\documentclass{IEEEtran}

\usepackage{graphicx}
\usepackage{acronym}
\usepackage{amsmath}
 \usepackage{url,doi} 

\usepackage{cite}
\usepackage{amssymb}

\usepackage{color}
\usepackage{balance}

\usepackage{bibunits}

\defaultbibliography{references}
\defaultbibliographystyle{IEEEtran}

\usepackage[ruled,linesnumbered,norelsize]{algorithm2e}

\makeatletter
\renewcommand{\@algocf@capt@plain}{above}
\makeatother
\usepackage[T1]{fontenc}

\newcommand{\measdim}{n_y} 
\newcommand{\statedim}{n_x} 
\newcommand{\mean}{\mu}
\newcommand{\statemean}{{\mu}}

\newcommand{\posteriormean}{\mean^+}

\newcommand{\weight}{w}

\newcommand{\meas}{y}

\newcommand{\predmeas}{\hat y}

\newcommand{\measfun}{h}

\newcommand{\statecov}{{\Pxx}} 
\newcommand{\posteriorcov}{{\statecov^+}} 

\newcommand{\noisecov}{R} 

\newcommand{\matr}[1]{\begin{bmatrix}#1\end{bmatrix}}

\newcommand{\state}{x}

\newcommand{\priorstate}{x_0}
\newcommand{\measnoise}{\varepsilon}

\newcommand{\jacobian}{J}

\newcommand{\constant}{b}

\newcommand{\covariance}{P}

\newcommand{\Pyy}{{\covariance_{\measfun{(\state)}\measfun{(\state)}}}}
\newcommand{\Pxx}{\covariance}
\newcommand{\Pxy}{{\covariance_{\state\measfun{(\state)}}}}
\newcommand{\Pyx}{{\covariance_{\measfun{(\state)\state}}}}

\newcommand{\priormean}{\mean^-}
\newcommand{\priorcov}{{\statecov^-}}

\newcommand{\nonlcov}{\Omega}

\newcommand{\innocov}{S}
\newcommand{\kalmangain}{K}

\newcommand{\nonlnoise}{{\varepsilon_\nonlcov}}

\newtheorem{Example}{Example}

\DeclareMathOperator{\N}{N}

\renewcommand{\d}{\mathrm{d}}

\SetCommentSty{mycommfont}

\begin{document}
\begin{bibunit}
\newacro{BinoGMF}{binomial Gaussian mixture filter}
\newacro{CDF}{cumulative distribution function}
\newacro{KF}{Kalman filter}
\newacro{KFE}{Kalman filter extension}

\newacro{CKF}{cubature Kalman filter}
\newacro{EKF}{extended Kalman filter}
\newacro{MC}{Monte Carlo}
\newacro{nEKF2}{numerical second order Extended Kalman filter}
\newacro{EKF2}{second order Extended Kalman filter}
\newacro{IPLF}{iterated posterior linearisation filter}
\newacro{GGF}{general Gaussian filter}
\newacro{GN}{Gauss-Newton}
\newacro{GMF}{Gaussian mixture filter}
\newacro{IKF}[IEKF]{iterated extended Kalman filter}
\newacro{IPLF}{iterated posterior linearization filter}

\newacro{LM}{Levenberg-Marquadt}

\newacro{DIPLF}{damped iterated posterior linearization filter}\newacro{KL}{Kullback-Leibler}
\newacro{KLD}{Kullback-Leibler divergence}
\newacro{KLPUKF}{Kullback-Leibler partitioned update Kalman filter}
\newacro{PUKF}{partitioned update Kalman filter}
\newacro{PDF}{probability density function}
\newacro{PF}{particle filter}
\newacro{RPLF}{regularized posterior linearization filter}
\newacro{DPLF}[DIPLF]{damped iterated posterior linearization filter}

\newacro{SLF}{statistically linearized filter}

\newacro{SLR}{statistical linear regression}

\newacro{S2KF}[S\textsuperscript{2}KF]{smart sampling Kalman filter}

\newacro{RUF}{recursive update filter}
\newacro{UKF}{unscented Kalman filter}
\newacro{MAP}{maximum a posteriori}

\title{Damped Posterior Linearization Filter}
\author{{Matti Raitoharju, Lennart Svensson, \'Angel F. Garc\'ia-Fern\'andez, and Robert Pich\'e}
\thanks{This article is carried out with support from a mobility grant from COST Action TU1302 (SAPPART).
 Financial support by the Academy of Finland under grant no. \#287792 (OpenKin) is hereby gratefully acknowledged.}
\thanks{M. Raitoharju is with the laboratory of Automation and Hydraulic Engineering, Tampere University of Technology and Department of Electrical Engineering and Automation, Aalto University, Finland. E-mail: {matti.raitoharju@tut.fi}. L. Svensson is with the department of Signals and Systems, Chalmers University of Technology, Sweden. Á. F. García-Fernández is with the Department of Electrical Engineering and Electronics, University of Liverpool, UK, Robert Pich\'e is with the laboratory of Automation and Hydraulic Engineering, Tampere University of Technology, Finland }
\thanks{This is an author accepted version of a letter published in IEEE Signal Processing Letters: \href{https://doi.org/10.1109/LSP.2018.2806304}{DOI:10.1109/LSP.2018.2806304} \copyright 2018 IEEE. }
}
\pagenumbering{arabic}

\maketitle

\begin{abstract}
In this letter, we propose an iterative Kalman type algorithm based on posterior linearization. The proposed algorithm uses a nested loop structure to optimize the mean of the estimate in the inner loop and update the covariance, which is a computationally more expensive operation, only in the outer loop. The optimization of the mean update is done using a damped algorithm to avoid divergence. Our simulations show that the proposed algorithm is more accurate than existing iterative Kalman filters.

\end{abstract}
\acresetall

\section{Introduction}
In Bayesian state estimation, a state that evolves stochastically in time is estimated from noisy measurements. In this letter, we concentrate on the measurement update stage of the Bayesian filter. In the measurement update, a prior (the dynamic model's state prediction) is updated using  information from a measurement. If the measurement model is linear and Gaussian the posterior density can be computed analytically using the Kalman filter \cite{kalman}, but for a general measurement model, the computation of the posterior density is intractable. One approximate approach is the \ac{GGF}, which represents the joint state and measurement distribution by a multivariate Gaussian using moment matching,  and then uses standard marginalization and conditioning formulas to compute the conditional state distribution given the measurement  to obtain the posterior estimate \cite{sarkka2013bayesian}. The \ac{GGF} moment matching is also intractable, but there exist several \acp{KFE} that are approximations of the \ac{GGF}. The standard way to apply \ac{GGF} is to do \ac{SLR}, or an approximation of the \ac{SLR}, in the prior; this can be interpreted as the optimal linearization with the given prior. The \ac{GGF} approach does not work well for severe nonlinearities or low measurement noise \cite{6584787} so there is interest in developing alternatives.

Iterative algorithms generate approximative solutions based on previous solutions with the goal of improving the approximation. The \ac{IKF} \cite[pp. 349-351]{jazwinski} produces a sequence of mean estimates by making first-order Taylor approximations of the measurement function. In \cite{Bell:1993}, it was shown that the \ac{IKF} measurement update is equivalent to the \ac{GN} algorithm for computing the \ac{MAP} estimate. \ac{IKF} performs well in situations where the true posterior is close to a Gaussian. Convergence of the \ac{GN} algorithm is not guaranteed. Better convergence is offered by damped (or descending) \ac{GN} algorithms \cite{Bjork}, which ensure that the cost function is nonincreasing at each iteration.
An alternative to damped \ac{GN} is the Levenberg-Marquadt algorithm \cite{Bellaire}.

There are many other iterated filters in the literature. The iterated sigma point Kalman filter \cite{Sibley-RSS-06}, at each iteration, uses a linearization
that mixes \ac{SLR} with respect to the prior and analytical linearization at the current \ac{MAP} estimate. The \ac{RUF} uses several updates of the prior with down-weighted Kalman gain \cite{UKFRUF}. Progressive Gaussian Filtering is a homotopy continuation method that updates the prior starting from an easily computable measurement likelihood that gradually evolves into the true measurement likelihood \cite{6290507,hanebeck2003progressive,7402897}. The \ac{KLPUKF} can be used when the measurement is multidimensional \cite{raitoharju2016}. It sequentially updates using nonlinearity-minimizing linear combinations of the measurements.

In this paper, we focus on the  \ac{IPLF}, which uses \ac{SLR} w.r.t. the posterior density. \ac{SLR}'s linearization is based on a larger area determined by a \ac{PDF} instead of a point as in the Taylor linearization, which improves the accuracy of the filter.  In~\cite{PLF}, it was shown that it is better for the linearization to be w.r.t. the posterior instead of the prior.  The basic idea in the \ac{IPLF} is to compute a posterior estimate using the prior, then use this estimate to compute a better linearization and posterior estimate, and so on. In this letter, we observe that the \ac{IPLF} does not always converge, and we propose a damped version of the algorithm with improved convergence properties.

We show in simulations how the proposed algorithm is less prone to filter divergence than the original \ac{IPLF}. Furthermore, we compare the posterior estimate accuracy with other iterative Kalman filter methods and show that the proposed algorithm is more accurate in the simulations.

\section{Background}
We consider measurements of the form \begin{equation}
\meas=\measfun(\state)+\measnoise, \label{equ:meas1}
\end{equation}
where $\meas$ is the $\measdim$ dimensional real measurement value, $\measfun(\cdot)$ is the measurement function, $\state$ is the $\statedim$ dimensional real random state vector and $\measnoise$ is zero mean Gaussian noise with covariance~$\noisecov$.
\subsection{\ac{GGF}}

The \ac{GGF} \cite{sarkka2013bayesian} for measurement  model \eqref{equ:meas1} uses the expected value of predicted measurement $\predmeas$, cross covariance of measurement and state $\Pyx$, and  measurement covariance $\Pyy$. These moments are
\begin{align}
	\predmeas &= \int \measfun(\state)p_{\N}(\state ; \priormean, \priorcov) \d \state \label{equ:i1}\\
	\Pyx &=\int  \left( \measfun(\state)  - \predmeas\right)\left( \state - \priormean\right)^Tp_{\N}(\state ; \priormean, \priorcov) \d \state \label{equ:i2} \\
	\Pyy\!&=\!\int\!\left( \measfun(\state)  - \predmeas\right)\left( \measfun(\state)  - \predmeas\right)^T\!p_{\N}(\state ; \priormean, \priorcov)  \d \state \label{equ:i3},
\end{align}	
where $p_{\N}(\state | \priormean, \priorcov)$ is the multivariate normal \ac{PDF} of the prior that has mean $\priormean$ and covariance $\priorcov$. The update equations for the mean and covariance of the posterior are then 
\begin{align}
	\posteriormean & = \priormean + \kalmangain(\meas - \predmeas) \\
	\posteriorcov & = \priorcov - \kalmangain\innocov\kalmangain^T, 
\end{align}
where
\begin{align}
	\innocov & =\Pyy + \noisecov  \\ 
	\kalmangain &= \Pxy \innocov^{-1} \label{equ:kgain}.
\end{align}

The integrals \eqref{equ:i1}-\eqref{equ:i3} do not have in general closed form solutions. Many \acp{KFE}, such as \ac{UKF} \cite{WANUKF} and \ac{CKF} \cite{cubature}, can be interpreted as numerical approximations of \ac{GGF}.

\subsection{\ac{IPLF}}
\label{sec:IPLF}
In \cite{PLF}, it is shown that a lower bound of the \ac{KLD} of a \ac{GGF} with respect to the true posterior is minimized if the moments \eqref{equ:i1}-\eqref{equ:i3} are computed using the posterior instead of the prior. This update cannot be computed directly with the \ac{GGF} equations \eqref{equ:i1}-\eqref{equ:kgain}, but instead a linear model that has correct moments is generated and is then applied to the prior. 

The \ac{IPLF} algorithm iteratively approximates the posterior using \ac{SLR}, creating a sequence of updated means $\statemean_i$  and covariance matrices $\statecov_i$ as follows. First, set $\statemean_1=\priormean$ and  $\statecov_1=\priorcov$. Then, at the $i$th iteration, compute the \ac{SLR} of $h(\cdot)$ with respect to ($\statemean_i$, $\statecov_i$) by first computing moments \eqref{equ:i1}-\eqref{equ:i3} and then defining the linear measurement that corresponds to those moments
\cite{PLF}:
\begin{equation}
	\hat h(x) = \jacobian_i \state + \constant_i + \nonlnoise_i + \measnoise , \label{equ:meas2}
\end{equation}
where
\begin{align}
\jacobian_i &= \Pxy_i^T \Pxx_i^{-1} \label{equ:jacobian}\\ 
\constant_i &=  \predmeas_i - \jacobian_i \statemean_{i}  \\
\nonlcov_i &=  \Pyy_i - \jacobian_i \Pxx_i \jacobian_i^T \label{equ:noncov} \\ 
\nonlnoise_i &\sim \N(0, \nonlcov_i),
\end{align}
$\nonlnoise_i$ is an independent noise whose covariance $\nonlcov_i$ is the covariance of the linearization error. For linear systems $\nonlcov_i$ is 0.

Finally, compute the posterior mean and covariance that correspond to the linearized measurement function \eqref{equ:meas2}:
\begin{align}
	\innocov &= \jacobian_i \priorcov \jacobian_i^T + \noisecov + \nonlcov_i \\ \
	\kalmangain_i & =  \priorcov \jacobian_i^T \innocov^{-1}\\
	\statemean_{i+1} & = \priormean + \kalmangain_i(\meas - \jacobian_i\priormean - \constant_i) \\
	\statecov_{i+1} & = \priorcov - \kalmangain_i\innocov_i\kalmangain_i^T. \label{equ:posteriorcov}
\end{align}
The obtained posterior is used for the next linearization. The process is repeated for a predetermined number of steps or until a convergence criterion is met, e.g.\ until the \ac{KLD} of two consecutive estimates is below a threshold.

However, the \ac{IPLF}, like the \ac{GN}, can diverge, as illustrated in the next example.

\begin{Example}
The measurement model is 
\begin{equation}
	\meas=\arctan(\state) + \varepsilon.
\end{equation}
State $\state$ has prior mean 2.75 and variance 1, the measurement value is $y=0$ and measurement noise variance is $R=10^{-4}$. The integrals \eqref{equ:i1}-\eqref{equ:i3} are computed using \ac{MC} integration with $10^5$ samples and also with \ac{IKF}.
\begin{table}
\centering
\caption{Mean estimates of first 6 iterations of \ac{IPLF} and \ac{IKF}}
\label{tbl:ex1}
\begin{tabular}{c|cccccc}
Iteration& 1 &2 &3 &4 &5 &6 \\ \hline
\ac{IPLF} &-2.51 &   5.28 & -31.75 &    17.52 &  -40.56 &   11.96 \\  
\ac{IKF} & -7.64 &   58.28 &   -1.77 &    2.60 &   -6.66 &   48.47  \\
\end{tabular}
\end{table} Table~\ref{tbl:ex1} gives the means in the first 6 iterations. Neither of the algorithms converge to the true posterior, which is close to 0, even after 50 iterations.
In the supplementary material, we provide another example, whose moments have closed form solutions, where the \ac{IPLF} does not converge.

\end{Example}

\subsection{\ac{IKF} as a \ac{GN} algorithm}
\ac{IKF} is similar to \ac{IPLF}, but it uses a first order Taylor series approximation of the measurement function instead of \ac{SLR}.   The \ac{IKF} iteration formulas are as in Section~\ref{sec:IPLF}, but with \eqref{equ:jacobian}-\eqref{equ:noncov} replaced by
\begin{align}
\jacobian_i &= \left.\frac{\d \measfun(x)}{\d x} \right|_{\statemean_{i-1}} \label{equ:IKF1}\\ 
\constant_i &= \measfun(\mean_{i-1}) - \jacobian_i \statemean_{i}  \\
\nonlcov_i &= 0. \label{equ:IKF3}
\end{align}
In \cite{Bell:1993}, it is shown that the \ac{IKF} for measurements of form \eqref{equ:meas1} is equal to the \ac{GN} algorithm \cite{Bjork} that minimizes the cost function
\begin{equation}
\begin{aligned}
	q(\statemean)= & \frac{1}{2}(\measfun(\statemean) - \meas)^T R^{-1}(\measfun(\statemean) - \meas) \\ &+ \frac{1}{2}(\statemean - \priormean)^T\statecov^{-1} (\statemean - \priormean) \label{equ:IKFcost}.
	\end{aligned}
\end{equation}

To decrease the possibility of the divergence of \ac{GN}, the damped (or descending) version of the \ac{GN} algorithm can be used \cite{Bjork}. A damped \ac{IKF} that uses line search \cite{Bjork} in minimization of \eqref{equ:IKFcost} is given in \cite{7266781}. The mean update is
\begin{equation}
\statemean_{i+1} = (1-\alpha)\statemean_{i}+ \alpha(\statemean_0 + K_{i}( y - \measfun(\statemean_{i} )))   \label{equ:meanup}
\end{equation}
with step length parameter $\alpha$ selected so that
$0\leq\alpha\leq1$ and  $q(\statemean_{i+1}) \leq q(\statemean_{i})$.
Finally, in the \ac{IKF}, the posterior covariance matrix is computed based on the linearization at the \ac{MAP} point after the last iteration  \cite{Bell:1993}.

\section{\ac{DPLF}}
The main weakness of the IPLF is that it sometimes diverges. In this section, we present a damped version of the IPLF, which has substantially better convergence properties. The key insight is that the update of the mean can be seen as an optimization problem and solved with a damped GN algorithm that decreases the cost function in every iteration. A detailed motivation of the algorithm is presented in Section~\ref{sec:derivation}, followed by a description of line search required to implement it. Further implementation details are discussed in the Supplementary material.

\subsection{Derivation of  \ac{DPLF}}
\label{sec:derivation}
Ideally, we want the posterior mean and covariance to be used in the \ac{SLR} to linearize the measurement function \cite{PLF}. 
Our approach to iteratively finding the mean and covariance for the linearization is to use a separate loop for the optimization of the mean. Thus,  the proposed filter computes the estimate using two nested loops. The inner loop uses a damped \ac{GN} algorithm to optimize the mean $\statemean_i$ while keeping  $\statecov_j $ and $\Omega_j$ fixed. The outer loop updates the covariances $\statecov_j $ and $\Omega_j$. Index $i$ increases in every iteration of the inner loop and $j$ in the outer loop.

\subsubsection{Inner loop}
The purpose of the inner loop is to find optimal $\statemean_i$ when $\statecov_j$ and $\nonlcov_j$ are fixed. At an optimum $\mu_i$ coincides with the mean used in the \ac{SLR} and we can write:
\begin{equation}
\begin{aligned}
	\statemean_i & = \priormean + \kalmangain(\statemean_i)(\meas - \jacobian(\statemean_i)\priormean - \constant(\statemean_i)),
\end{aligned}
\end{equation}
where we used $(\statemean_i)$ to identify vectors and matrices that depend on the mean. In the supplementary material, we show that the above equation is fulfilled when  
\begin{equation}
\begin{aligned}
	q(\statemean_i)= & \frac{1}{2}( \hat{y}(\statemean_i)- \meas)^T (R+\Omega_j)^{-1}( \hat{y}(\statemean_i) - \meas) \\& + \frac{1}{2}(\statemean_i - \priormean)^T\left( \priorcov \right)^{-1} (\statemean_i - \priormean) \label{equ:target}
\end{aligned}
\end{equation}
achieves its minimum. Therefore, we use it as the cost function for the inner loop. 
The optimization problem is analogous to the optimization of the cost function of the \ac{IKF} \eqref{equ:IKFcost}, but using $\hat{y}$ instead of $\measfun(\statemean)$. Because the step is taken in a descent direction of the target function we can use a scaling factor $\alpha$ ($0<\alpha \leq 1$) in the update  as in~\eqref{equ:meanup}:
\begin{equation}
	\mu_{i+1}  = (1-\alpha) \mu_i + \alpha \left( \mu_0 + K_i(\meas - J_i \mu_0 - b_i) \right) \label{equ:Kalmanalpha} 
\end{equation}
and select $\alpha$ using line search so that the value of our target function \eqref{equ:target} decreases. Thus, the next mean is a weighted sum of previous mean and the mean provided by \ac{IPLF}.  This makes the algorithm  locally convergent on almost all nonlinear least squares problems, provided that the line search is carried out appropriately. In fact, it is usually globally convergent \cite[Chapter 9.2.1]{Bjork}.  To facilitate the stopping of the inner loop, we propose that the algorithm should repeat the inner loop as long as \eqref{equ:target} reduces significantly at each iteration i.e.\ while
\begin{equation}
	q(\statemean_{i+1},\statecov_j) < \beta q(\statemean_i,\statecov_j)
\end{equation}
with, e.g.,\ $\beta=0.9$.

\subsubsection{Outer loop}
The outer loop updates the covariances $\statecov_j$ and $\nonlcov_j$ for the next iteration of the inner loop. To define a stopping condition for the outer loop, we first exponentiate function \eqref{equ:target} and normalize to get a product of normal distributions:
\begin{equation}
\begin{aligned}
	& \frac{1}{\sqrt{(2\pi)^{\measdim} |R+\Omega_j|(2\pi)^{\statedim}|\priorcov|}}e^{-q(x)} \\ = & 	p_{\N} ( \predmeas_i | \meas, R+\Omega_j) p_{\N}( \statemean_i | \priormean, \priorcov) \label{equ:ML}.
	\end{aligned}
\end{equation}
If we treat $P_i$ and  $\Omega_i$ as constants, \eqref{equ:target} achieves its minimum when \eqref{equ:ML} achieves its maximum. We thus continue the outer loop as long as \eqref{equ:ML} is increasing significantly. 
Compared to using \eqref{equ:target}, the factor containing  $\Omega_i$ in \eqref{equ:ML} causes the stopping condition based on \eqref{equ:ML} to favor estimates with smaller nonlinearity.

\subsubsection{Full algorithm}

The algorithm is summarized in Algorithm~\ref{algo:gnplf}.
\begin{algorithm}[tb]
\small
\caption{\ac{DPLF}}
\label{algo:gnplf}
	$i\leftarrow0$, $j\leftarrow0$,  $\mu_0 \leftarrow \priormean$, $\statecov_0 \leftarrow \priorcov$ \\
	Compute  $\predmeas_i$, $\jacobian_i$, and $\Omega_j$ at prior using \eqref{equ:i1}, \eqref{equ:jacobian}, and \eqref{equ:noncov}\\
	\While{$p_{\N} ( \predmeas_i | \meas, R+\Omega_j) p_{\N}( \mean_i | \priormean, \priorcov)$  increases significantly}
	{
	\While{$q(\mu_i,P_j) < \beta q(\mu_{i-1},P_j)$ }
	{
	$i \leftarrow i +1$ \\
	Find $0<\alpha_i \leq1$ for \eqref{equ:Kalmanalpha} so that \eqref{equ:target} becomes smaller using line search (Section~\ref{sec:LS})\\
	 Compute $\statemean_i$ using \eqref{equ:Kalmanalpha} \\ 
	 Compute $\jacobian_i$ using \eqref{equ:i2} and \eqref{equ:jacobian} \\
	 Compute $\predmeas_i$ using \eqref{equ:i1}  
	 }
 		 $j \leftarrow j  +1$ \\

	Compute $\statecov_j$ using \eqref{equ:posteriorcov} \\
	Compute $\nonlcov_j$ using \eqref{equ:noncov}

	}
	Use the mean and covariance from the second last round of the outer loop, i.e.\ those that had highest $p_{\N} ( \predmeas_i | \meas, R+\Omega_j) p_{\N}( \mean_i | \priormean, \priorcov)$, as the posterior estimate
\end{algorithm}
The inner loop does the optimization of the target function \eqref{equ:target} while keeping the state covariance and nonlinearity covariance constant. The direction of change in a \ac{GN} algorithm is a descent direction of the cost function \cite[Chapter 9.2.1]{Bjork}. Thus, with small enough $\alpha$ the cost function value should decrease, unless the estimate is in a local extremum or a saddle point.  By using $\alpha=1$ and exiting the inner loop after a single iteration, Algorithm~\ref{algo:gnplf} becomes equivalent to an \ac{IPLF} algorithm.

The outer loop changes the optimization problem of the inner loop as $\nonlcov$ changes and the covariance, for which the expected measurement value $\predmeas$ is computed, changes. 

\subsection{Line search}
\label{sec:LS}
The value of $\alpha$ can be computed with different line search techniques \cite[Chapter 9.2.1]{Bjork}. Backtracking line search, also known as the Armijo-Goldstein step length principle, is one of the simplest and most widely used options, and is done as follows.

First, $\alpha$ is set to 1, and a value of $\statemean$  is found using \eqref{equ:meanup}. If this value does not decrease \eqref{equ:target}, w.r.t. the previous value of $\statemean$, then $\alpha$ is multiplied by a constant factor $0<\tau<1$ and the procedure is repeated. To reduce the computational load one can terminate the inner loop when $\alpha$ is smaller than a predetermined value (we use $2^{-4}$) as then the change is the estimate is negligible.

\section{Simulations}
\label{sec:simulations}
In our first test, we evaluate the situation in Example~1.
We compare posterior estimates obtained with estimation algorithms that are based on approximating integrals \eqref{equ:i1}-\eqref{equ:i3}. We use \ac{MC} integration with $10^5$ samples, \ac{EKF}, \ac{UKF} with parameters from \cite{WANUKF}, and third-degree \ac{CKF} \cite{cubature} to compute the moment approximations. The algorithms used are \ac{GGF}, \ac{RUF} \cite{UKFRUF}, \ac{IPLF} \cite{PLF}, and \ac{DPLF}.  The \ac{GGF} with different moment computation algorithms is equivalent to the algorithms from which the moment computation was taken. \ac{RUF} applies the measurement using down-weighted Kalman gains and uses 10 iterations in our tests. When \ac{EKF} is used to compute the moments \ac{GGF} is \ac{EKF}, \ac{RUF} is equivalent to the algorithm in  \cite{RUKF}, \ac{DPLF} is \ac{IKF}, and \ac{DPLF} is damped \ac{IKF}. We use $\tau=0.5$ for reducing the step length.  For the threshold of inner loop, we use $\beta=0.9$ and we exit the outer loop if  $0.999 p_{\N} ( \predmeas_i | \meas, R+\Omega_j)$ is not larger than the corresponding likelihood of the previous iteration. A more thorough discussion on the selection of parameters is given in the supplementary material.
\begin{table}
\centering
\caption{\acp{KLD} of estimates produced by different filters in $\arctan$ test}
\label{tbl:klds1}
\begin{tabular}{c|ccccc}
& GGF  & \ac{RUF} & IPLF & \ac{DPLF}  \\ \hline
\ac{MC} &   15.9  &    0.04 &   88.69 &    $3\cdot 10^{-6}$ \\
\ac{EKF} &4009.10&    0.01 &   65.12 &   $10^{-6}$ \\
\ac{CKF} &3370.78&    0.01 &   64.39 &    $10^{-6}$   \\
\ac{UKF} &   92.55  &    0.01 &    $ 10^{-6}$ &   $ 10^{-6}$\\
 \end{tabular}
\end{table}

Results are given in Table~\ref{tbl:klds1}. \ac{DPLF} gives good posterior estimates with all moment computation methods.  \ac{IPLF} converges towards an estimate within 50 iterations only when using \ac{UKF} for moment computation. Though the \ac{UKF} and \ac{CKF} usually have similar accuracy, in this case, the \ac{IPLF} did not converge when using the \ac{CKF}; changing the prior mean a little may make it to converge with different moment computation algorithms. \ac{RUF} has slightly worse estimates than \ac{DPLF} and \ac{GGF} does not provide good estimates with any of the tested moment computation methods. 

In our second test, we computed estimates for two dimensional positioning using three range measurements in a single measurement update. The prior has zero mean with covariance $I$. Ranges were computed to beacons located at $\matr{-1&0}^T,\matr{0&1}^T$ and $\matr{1&-2}^T$. Measurement noise covariance was $I$ and measurements were generated by sampling a true location from the prior and then generating corresponding measurements. In this test, we used the same filters as in the previous test and \ac{KLPUKF} \cite{raitoharju2016}. \ac{KLPUKF} uses a linear transformation to decorrelate measurement elements in such a way that the nonlinearity of a measurement element is minimized and applies measurement elements sequentially. For the first test with only a scalar measurement, \ac{KLPUKF} would have been identical to \ac{GGF}.

The test is repeated 1000 times using different true location and measurement values. We analyze only state after one update step, but a more accurate update step can be expected to lead to improved stability in the filtering recursion. The mean \acp{KLD}, which are computed numerically using a dense grid, are presented in Table~\ref{tbl:klds2}. The iterative algorithms that use \ac{EKF} moment computation, except \ac{RUF}, produce worse estimates than \ac{EKF} itself. This is probably because the posterior may be multimodal in this test. With other moment computation methods the \ac{IPLF} is worse than \ac{GGF}, but \ac{DPLF} is better than either of those. In this test, \ac{KLPUKF} and \ac{DPLF} were the most accurate algorithms, except when tested with \ac{EKF}.  The residuals of the mean estimates compared to true locations were similar for all algorithms, thus the largest improvement of using \ac{DPLF} comes from better posterior covariance estimation. In this test, the computational time of \ac{DPLF} was 8 times higher than the computational cost of \ac{KLPUKF}. Our implementations of other algorithms had runtimes between these, \ac{IPLF} taking similar time as \ac{DPLF}.  The exact comparison of computational complexity of algorithms evaluated in this letter cannot be done, because algorithms have different conditions for loops. These conditions depend on the estimation problem and computational complexity cannot be stated as a function of dimensionalities of the state and the measurement vector.

\begin{table}
\centering
\caption{Mean \acp{KLD} of estimates in range test}
\label{tbl:klds2}
\begin{tabular}{c|cccccc}
& GGF & KLPUKF & RUF & IPLF & \ac{DPLF} \\ \hline
\ac{MC} &0.25 & 0.17 & 0.21 & 0.26 & 0.17  \\
\ac{EKF} &0.48 & 0.60 & 0.48 & 0.55 & 0.55   \\
\ac{CKF} &0.28 & 0.22 & 0.29 & 0.38 & 0.23  \\
\ac{UKF} & 0.35 & 0.34 & 0.33 & 0.37 & 0.26  
\end{tabular}
\end{table}

\section{Conclusions and future work}

We have shown that by fixing the covariance matrices of the \ac{IPLF} algorithm, the optimal posterior mean in the posterior linearization approach is the solution of an optimization problem.  We solved this optimization problem using a damped \ac{GN} algorithm. Then we updated the covariance matrices $\statecov_j$  and $\nonlcov_j$ and repeated the mean optimization. Compared to \ac{IPLF}, this algorithm is less prone to diverge. The proposed algorithm was more accurate than other tested algorithms in our simulations.

One topic for future research is to see how posterior linearization could be done with other optimization algorithms than \ac{GN}, such as Levenberg-Marquadt. Another interesting research topic would be extending the algorithm  to work with non-additive noise. We have done some preliminary simulations in time series, where we found that the proposed algorithm is more accurate than other iterative algorithms,  but all iterative algorithms have a risk of converging to a wrong local solution. Investigating this is one future topic.
\balance
\putbib
\end{bibunit}
\clearpage
\begin{bibunit}
\pagenumbering{Alph}
\onecolumn
\setcounter{equation}{0}
\renewcommand{\theequation}{S-\arabic{equation}}
\renewcommand{\thefigure}{S-\arabic{figure}}
\section*{Supplementary Material to Damped Posterior Linearization Filter} 
\subsection{Analytical example of a situation where the \ac{IPLF} does not converge}

\begin{figure}[b!]
\centering
\includegraphics[width=10cm]{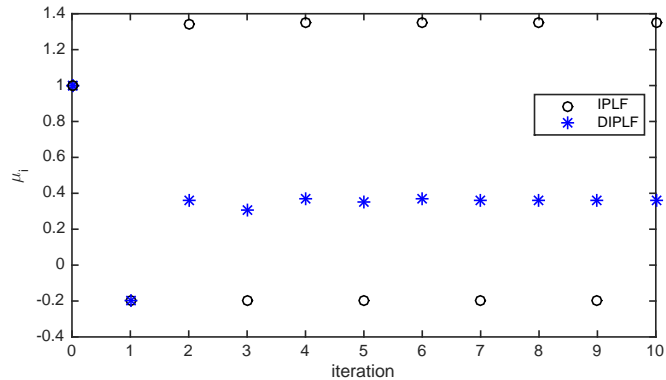}
\caption{Evolution of estimates of \ac{IPLF} and \ac{DIPLF}}
\label{fig:convergence}
\end{figure}

\begin{figure}[b!]
\centering
\includegraphics[width=11cm]{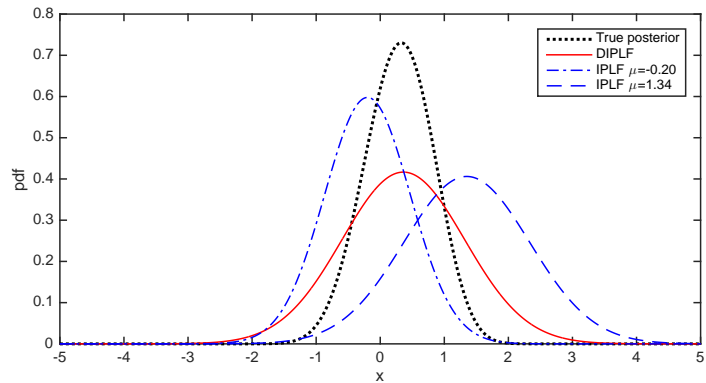}
\caption{True posterior compared to the oscillating estimates of \ac{IPLF} and the estimate provided by \ac{DIPLF}}
\label{fig:posterior}
\end{figure}

Consider a one dimensional state with prior $\priorstate \sim \N(\priormean,\priorcov)$ with $\priormean=1$ and $\priorcov=1$ and a measurement model 
\begin{equation}
	h(x)=y^2 + \measnoise, \label{equ:spoly}
\end{equation}
where $\measnoise \sim  \N(0,R)$ with $R=2^2$.  Because the measurement function is a second order polynomial, the second order \ac{EKF} \cite{sarkka2013bayesian} produces closed form \ac{SLR}. For \eqref{equ:spoly} the moments (2)-(4) in paper are
\begin{align}
	\hat{y} & = \statemean_i  + \statecov_j \\
	\Pyx &= 2\statemean_i\statecov_j \\
	\Pyy &= 2\statecov_j^2 .
\end{align}
Consider a measurement with value $y=-4$. After the first iteration, the mean of \ac{IPLF} is $-0.2$ and, after the third iteration, the mean is $1.34$, thus the mean has moved to other side of the initial mean after the first iteration. When the iterations are continued, the estimates' means oscillate between $-0.20$ and $1.35$. On the other hand, the \ac{DIPLF} has the same estimate after the first iteration, but after that it converges to a single estimate with mean $0.36$. Figure~\ref{fig:convergence} shows how the mean estimates converge during the iterations in \ac{IPLF} and \ac{DIPLF}.
Figure~\ref{fig:posterior} shows the \acp{PDF} of the true posterior, the oscillating estimates of the \ac{IPLF} and the estimate made by \ac{DIPLF}. Figure shows how the mean and mode of the estimate provided by the \ac{DIPLF} are close to the mean and mode of the true posterior. However, due to the nonlinearity of the measurement function, \ac{DIPLF} overestimates the covariance.

\subsection{Derivation of the cost function}
Here we derive the cost function for updating the inner loop. First, we modify the Kalman update similarly as in \cite{Bell:1993} to get a \ac{GN} update from \ac{IKF} update, but taking  the $\Omega_j$ term of \ac{IPLF} into account:
\begin{equation}
\begin{aligned}
	\statemean_i  = &\priormean + \kalmangain(\statemean_i)(\meas - \jacobian(\statemean_i)\priormean - \constant(\statemean_i)) \\
     = &\priormean  + (\jacobian(\statemean_i)^T(R+\Omega_j)^{-1}\jacobian(\statemean_i) + \priorcov ^{-1})^{-1}\jacobian(\statemean_i)^T(R+\Omega_j)^{-1} (y - \jacobian(\statemean_i) \priormean - \hat y (\statemean_i)  + \jacobian(\statemean_i)\mu_i) \\
	 = & (\jacobian(\statemean_i)^T(R+\Omega_j)^{-1}\jacobian(\statemean_i) + \priorcov ^{-1})^{-1}(\jacobian(\statemean_i)^T(R+\Omega_j)^{-1}( y - \hat y (\statemean_i)   + \jacobian(\statemean_i)\mu_i ) + \priorcov ^{-1}\priormean ) \\
	 =& (\jacobian(\statemean_i)^T(R+\Omega_j)^{-1}\jacobian(\statemean_i) + \priorcov ^{-1})^{-1}(\jacobian(\statemean_i)^T(R+\Omega_j)^{-1}( y - \hat y (\statemean_i)   + \jacobian(\statemean_i)\mu_i ) +\priorcov ^{-1}\mu_i -\priorcov ^{-1}\mu_i + \priorcov ^{-1}\priormean )\\
	 = &\mu_i + (\jacobian(\statemean_i)^T(R+\Omega_j)^{-1}\jacobian(\statemean_i) + \left(\priorcov\right)^{-1})^{-1}(\jacobian(\statemean_i)^T(R+\Omega_j)^{-1}( y - \hat y (\statemean_i)  ) + \priorcov ^{-1}(\priormean-\mu_i)) \label{equ:GNup}.
	 \end{aligned}
\end{equation}
From the last row we can see that the \ac{IPLF} update for the mean with fixed $\Omega_j$ has a stationary point when
\begin{equation}
\jacobian(\statemean_i)^T(R+\Omega_j)^{-1}(\hat y (\statemean_i)- y) + \priorcov ^{-1}(\mu_i-\priormean) = 0. \label{equ:stable}
\end{equation}
Differentiating the expected value of a measurement
\begin{align}
	\predmeas(\statemean_i) &= \int \measfun(\state)p_{\N}(\state ; \statemean_i, \statecov_j) \d \state
\end{align}
 with respect to the mean~$\mu_i$ gives
\begin{equation}
\begin{aligned}
	 &\frac{\d \int \measfun(\state)p(\state ; \mu_i, \statecov_j) \d \state }{\d \mu_i} \\
	 =& \int \measfun(\state)\frac{\d p(\state ; \mu_i, \statecov_j)  }{\d \mu_i} \d \state 
	\\=& \int \measfun(\state)( x - \mu_i )^T p(\state ; \mu_i, \statecov_j)  \d \state \statecov_j^{-1} .
\end{aligned}
\end{equation}
	Noting that $ \int \hat y(\statemean_i)( x - \mu_i )^T p(\state ; \mu_i, \statecov_j)\d \state = 0 $ we can write
\begin{equation}
\begin{aligned}
	& \int (\measfun(\state)-\hat y(\statemean_i))( x - \mu_i )^T p(\state ; \mu_i, \statecov_j)  \d \state \statecov_j^{-1} = \Pyx \statecov_j^{-1}.
\end{aligned}
\end{equation}
From this we see that the Jacobian computed from \ac{SLR}
\begin{equation}
\jacobian(\statemean_i) = \Pxy_i^T \Pxx_j^{-1} \label{equ:jacobian2}\\ 
\end{equation} is the Jacobian of the expected value of the measurement i.e.\ $\jacobian(\statemean_i) = \frac{\d \hat y (\statemean_i) }{\d \mu_i}$. Using this we see that the left hand side of \eqref{equ:stable} is the derivative of the function:
\begin{equation}
\begin{aligned}
	q(\statemean_i)= & \frac{1}{2}( \hat{y}(\statemean_i) - \meas)^T (R+\Omega_j)^{-1}( \hat{y}(\statemean_i) - \meas) + \frac{1}{2}(\statemean_i - \priormean)^T\left(\priorcov\right)^{-1} (\statemean_i - \priormean) \label{equ:targetproof},
\end{aligned}
\end{equation}
which has a local extremum when  \eqref{equ:stable} holds. Thus, \eqref{equ:targetproof}  can be used as the cost function of the inner loop. 

\subsection{Implementation considerations}
\label{sec:implementation}

\subsubsection{Implementation with sigma-point filters}
Integrals (2)-(4) do not have closed form solutions in general. Sigma-point filtering is a common framework for approximating these integrals \cite{sarkka2013bayesian}. The formulas are
\begin{align}
	\predmeas & \approx \sum_{k=1}^{m} \weight_k \measfun(\mean + \Delta_k ) \\
	\Pyx & \approx \sum_{k=1}^{m} \weight_k \left(\measfun(\mean + \Delta_k )- \predmeas \right)\Delta_k^T \\
	\Pyy & \approx \sum_{k=1}^{m} \weight_k \left(\measfun(\mean + \Delta_k )- \predmeas \right)\left(\measfun(\mean + \Delta_k )- \predmeas \right)^T,
\end{align}
where $\weight_k$ is the weight of $k$th sigma-point, $m$ is the number of sigma points, and $\Delta_k$ is the translation of the $k$th sigma-point from the state mean.  The translation vectors $\Delta_k$ depend on the state dimension and the state covariance matrix $\statecov$ and their computation involves taking a matrix square root (e.g.\ Cholesky decomposition) of $\statecov$. In Algorithm~1 the state covariance matrix is constant in the inner loop and so the $\Delta_k$ vectors need to be recomputed only in the outer loop.

\subsubsection{Efficient computation of the Jacobian}
The updated state covariance (17) can be written in information form using the matrix inversion lemma:
\begin{equation}
	\statecov_{j+1} = \left( \jacobian_i^T (R+\Omega_j)^{-1}\jacobian_i + \left(\priorcov\right)^{-1}\right)^{-1}.
\end{equation}
Using this in (10) we get
\begin{equation}
	\jacobian_{i+1} = \Pyx \left( \jacobian_i^T (R+\Omega_j)^{-1}\jacobian_i + \left(\priorcov\right)^{-1}\right). \label{equ:opti}
\end{equation}
Since $\priorcov$ does not change between iterations,  $\left(\priorcov\right)^{-1}$ needs to be computed only once for the inner loop and  matrices that have to be inverted at every iteration of the inner loop have dimension $\measdim \times \measdim$. Thus, when $\measdim \ll \statedim$ \eqref{equ:opti} is faster than (10).

\subsubsection{Selection of parameters $\tau$ and $\beta$}

\begin{figure}
\includegraphics[width=\textwidth]{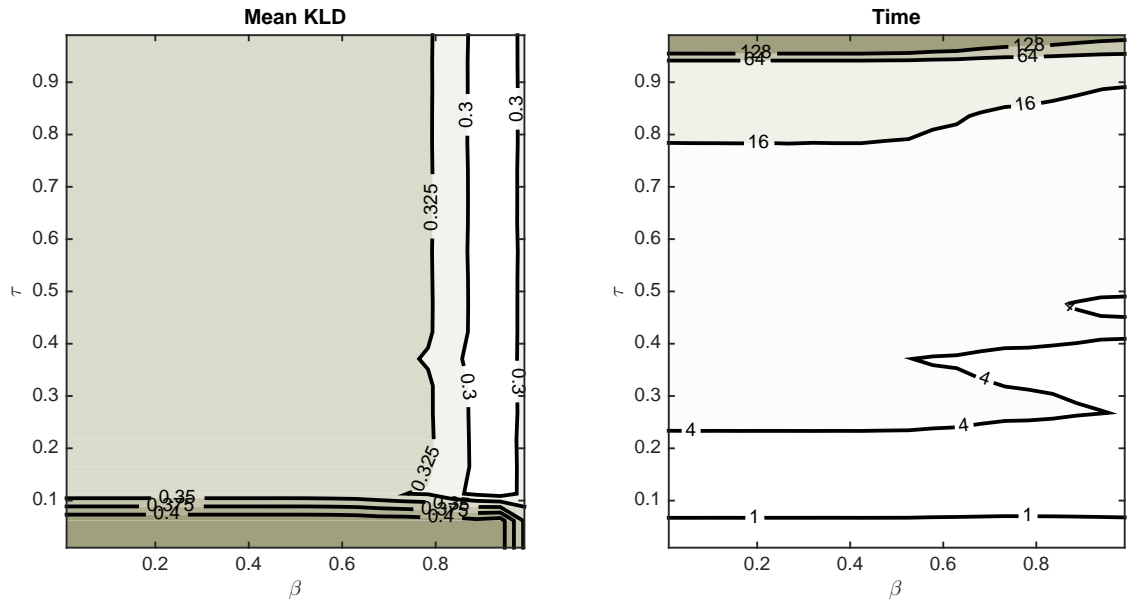}
\caption{Effect of parameters $\tau$ and $\beta$ on the mean \ac{KLD} and computing time}
\label{fig:runtime}

\end{figure}

The parameters $\tau$ and $\beta$ govern the behavior of inner loop of the \ac{DIPLF}. Setting a low $\beta$ makes the algorithm exit faster from the inner loop. Using a large $\tau$ changes the value of $\alpha$ less and makes the algorithm slower, while setting a low $\tau$ reduces the step length fast. Figure~\ref{fig:runtime} shows the effect of the parameters $\tau$ and $\beta$ on the \ac{KLD} and on the computing time in our second test in Section~IV. Figure~\ref{fig:runtime} shows that using $\tau < 0.12 $ reduces the accuracy, unless $\beta$ is close to 1. When $\tau$ is larger than $0.12$ there is not much difference in the accuracy, but the runtime increases fast. The parameter $\beta$  has a smaller effect on the runtime than $\tau$ and for the accuracy the optimal value is near $0.9$. Naturally, the given values are valid only for the example in this paper.

\putbib
\end{bibunit}

\end{document}